# Fast Approximate Polynomial Multipoint Evaluation and Applications

*Alexander Kobel*[1–3] *and Michael Sagraloff*[1]

✉ alexander.kobel@mpi-inf.mpg.de • michael.sagraloff@mpi-inf.mpg.de
[1] Max-Planck-Institut für Informatik • Campus E1 4, 66123 Saarbrücken, Germany
[2] International Max Planck Research School for Computer Science • 66123 Saarbrücken, Germany
[3] Universität des Saarlandes • 66123 Saarbrücken, Germany

**Abstract**

It is well known that, using fast algorithms for polynomial multiplication and division, evaluation of a polynomial $F \in \mathbb{C}[x]$ of degree $n$ at $n$ complex-valued points can be done with $\tilde{O}(n)$ exact field operations in $\mathbb{C}$, where $\tilde{O}(\cdot)$ means that we omit polylogarithmic factors. We complement this result by an analysis of *approximate multipoint evaluation* of $F$ to a precision of $L$ bits after the binary point and prove a bit complexity of $\tilde{O}(n(L + \tau + n\Gamma))$, where $2^\tau$ and $2^\Gamma$, with $\tau, \Gamma \in \mathbb{N}_{\geq 1}$, are bounds on the magnitude of the coefficients of $F$ and the evaluation points, respectively. In particular, in the important case where the precision demand dominates the other input parameters, the complexity is soft-linear in $n$ and $L$.

Our result on approximate multipoint evaluation has some interesting consequences on the bit complexity of further approximation algorithms which all use polynomial evaluation as a key subroutine. Of these applications, we discuss in detail an algorithm for polynomial interpolation and for computing a Taylor shift of a polynomial. Furthermore, our result can be used to derive improved complexity bounds for algorithms to refine isolating intervals for the real roots of a polynomial. For all of the latter algorithms, we derive near-optimal running times.

**Keywords**

approximate arithmetic
fast arithmetic
multipoint evaluation
certified computation
polynomial division
polynomial interpolation
Taylor shift
root refinement

We thank Victor Pan for making us aware of the extensive article "Partial Fraction Decomposition in $\mathbb{C}(z)$ and Simultaneous Newton Iteration for Factorization in $\mathbb{C}[z]$" by Peter Kirrinnis [11]. In a slightly different context, coined as "modular representation in $\mathbb{C}(z)$" [11, Algorithm 3.5 and Theorem 3.9], he shows accuracy and bit complexity estimates for approximate multipoint evaluation that are similar to our contribution. Unfortunately, this part of Kirrinnis' work seems to be almost unknown in the community. We assume that this is due to the fact that the theorem is not prominently highlighted, but only serves as an auxiliary lemma, with its condensed proof put in the appendix. Furthermore, Kirrinnis' article is virtually impossible to find in searches with the keywords "(approximate) multipoint evaluation," which do not appear at all in the publication.

Although this part of our contribution eventually turned out to be already known, we believe that the article at hand is a useful resource: it provides an in-depth and mostly self-contained explanation of the subject, and the results are stated in a more general form directly suitable for more applications (e.g., non-monic polynomials). Furthermore, we discuss several applications of approximate multipoint evaluation. We also want to mention our subsequent publication [12], where we use fast approximate multipoint evaluation to derive record bounds for several important problems in the study of planar algebraic curves, and analyze a different and conceptually simpler approach to the polynomial division subroutine that leads to a comparable bit complexity.

# 1. Introduction

We study the problem of approximately evaluating a polynomial $F \in \mathbb{C}[x]$ of degree $n$ at $n$ points $x_1, \ldots, x_n \in \mathbb{C}$. More precisely, we assume the existence of an oracle which provides arbitrarily good approximations of the polynomial's coefficients as well as of the points $x_i$ for free; that is, querying the oracle is as expensive as reading the approximation. Under this assumption, we aim to compute approximations $\tilde{y}_i$ of $y_i := F(x_i)$ such that $|y_i - \tilde{y}_i| \leq 2^{-L}$ for all $i = 1, \ldots, n$, where $L \in \mathbb{N}$ is a given non-negative integer. In what follows, let $2^{\tau}$ and $2^{\Gamma}$ with $\tau, \Gamma \in \mathbb{N}_{\geq 1}$ be upper bounds for the absolute values of the coefficients of $F$ and the points $x_i$, respectively.

When considering a sequential approach, where each $\tilde{y}_i \approx F(x_i)$ is computed independently by using Horner's Scheme and approximate but certified interval arithmetic [9], we need $O(n)$ arithmetic operations with a precision of $O(L + \tau + n\Gamma)$ for each of the points $x_i$. Thus, the total cost for all evaluations is bounded by $\tilde{O}(n^2(L + \tau + n\Gamma))$ bit operations.[1]

In this paper, we show that using an approximate variant of the classical fast multipoint evaluation scheme [14, 7], we can improve upon the latter bound by a factor of $n$ to achieve $\tilde{O}(n(L + \tau + n\Gamma))$ bit operations. The classical fast multipoint evaluation algorithm reduces polynomial evaluation at $n$ points to successive polynomial multiplications and divisions which are all balanced with respect to degree. It is a well known fact that, for exactly computing all values $y_i$, it uses only $O(n \log^2 n)$ *exact field operations* in $\mathbb{C}$ compared to $O(n^2)$ field operations if all evaluations are carried out independently; see Section 2 for a short review. This method has mostly been studied for low precisions, in particular for its performance with machine floating point numbers; see, e.g., [2, Section 2] or the extensive discussion in [13]. It is widely considered to be numerically unstable, mainly due to the need of polynomial divisions, and the precision demand for the sequential evaluations based on Horner's scheme does not directly carry over.

In previous work (e.g., [15, 16]), more involved algorithms for fast approximate multipoint evaluation have been introduced that allow to decrease the total number of (infinite precision) arithmetic operations from $O(n \log^2 n)$ to $O(n \log n)$ (if $n$ dominates all other input parameters). The authors mainly focus on the *arithmetic complexity* of their algorithms, and thus no bound on the *bit complexity* is given. For the special case where the points $x_i$ are the roots of unity, the problem can be solved with $\tilde{O}(n(\tau + L))$ bit operations by carrying out the fast Fourier transform with approximate arithmetic [20, Theorem 8.3]. However, for general points, we are not aware of any bit complexity result which considerably improves upon the bound $\tilde{O}(n^2(L + \tau + n\Gamma))$ that one directly obtains from carrying out all evaluations independently.

The main contribution of this paper is to show that the previously claimed issue of numerical instability within the classical fast multipoint evaluation scheme can be resolved. The crux of our approach is as follows: First, we exploit the fact that all divisors in the considered polynomial divisions are monic polynomials $g_{i,j}(x) = (x - x_{(j-1)\cdot 2^i+1}) \cdots (x - x_{j\cdot 2^i})$, which allow a stable division, at least if the precision $L$ dominates the values $n$ and $\Gamma$; see Corollary 8 for a precise statement. Second, we consider a numerical division algorithm from Schönhage [18] which yields an output precision of $L$ bits after the binary point if the algorithm runs with an *internal* precision of $\approx 2 \cdot L$. The crucial misconception in previous research is the conclusion that the input of the division algorithm must be available in a precision of $2 \cdot L$ bits as well. Hence, in the multipoint evaluation algorithm where we have to consider $\log n$ successive divisions, an accuracy of $n \cdot L$ seems to be necessary in the initial stages to guarantee $L$ meaningful bits at the end. However, we show that the precision requirement only holds for *intermediate* calculations, and that it suffices to consider only $L$-bit approximations of the input; see Section 2.2. Thus, the propagated error throughout the multipoint evaluation scheme stays within $\approx 2^{-L}$ compared to $\approx 2^{-L/n}$, effectively lowering the precision demand and, consequently, the bit complexity of the entire algorithm by a factor of $n$ upon the previously known bounds.

---
[1] $\tilde{O}(\cdot)$ means that polylogarithmic factors are ignored.

Our result immediately improves the bit complexity of many other approximation algorithms which heavily use polynomial evaluation as a key subroutine. One example and an important application of multipoint evaluation is *polynomial interpolation*. For given points $x_1, \ldots, x_n \in \mathbb{C}$ and corresponding interpolation values $v_1, \ldots, v_n$, there exists a unique polynomial $F \in \mathbb{C}[x]$ of degree less than $n$ such that $F(x_i) = v_i$ for all $i$. Based on our approach for fast multipoint evaluation, we prove that computing an $L$-bit approximation $\tilde{F}$ of $F$ (i.e., $\|\tilde{F} - F\|_1 \leq 2^{-L}$) uses only $\tilde{O}(nL)$ bit operations (for $L$ dominating $n$ and the bitsizes of the $x_i$'s and $v_i$'s). Our more general complexity bound as stated in Section 3.1 also involves the absolute values of the points $x_i$ and the values $v_i$ as well as the geometric location of the $x_i$'s.

Furthermore, we combine fast approximate multipoint evaluation and approximate interpolation in order to derive an alternative method to [18, Theorem 8.4] for computing an $L$-bit approximation of a *Taylor shift* of a polynomial $F$ (i.e., the polynomial $F_m(x) := F(m + x)$ for some $m \in \mathbb{C}$) with $\tilde{O}(nL)$ bit operations (again, for $L$ dominating). The details are given in Section 3.2.

Finally, approximate multipoint evaluation facilitates near-optimal algorithms for the simultaneous *refinement of real root approximations* for polynomials. More precisely, for a variant of the Quadratic Interval Refinement method [1, 9, 10] which uses approximate computation, we can directly improve the bit complexity for computing $L$-bit approximations of all real roots of $F \in \mathbb{R}[x]$ from $\tilde{O}(n^2 L)$ to $\tilde{O}(nL)$, with $n := \deg F$ (if $L$ dominates parameters that only depend on $F$ such as the separations of its roots). This improvement mainly stems from the fact that, instead of considering the refinements of each of the isolating intervals independently, we can merge all evaluations of $F$ in a certain precision in a single multipoint evaluation, with only logarithmic overhead compared to a single evaluation. A brief overview of this method is given in the Appendix A.2; for a detailed description, we refer to [10].[2]

## 2. Approximate Polynomial Multipoint Evaluation

Given a polynomial $F(x) = \sum_{i=0}^{n} F_i x^i \in \mathbb{C}[x]$ of degree $n$, complex points $x_1, \ldots, x_n \in \mathbb{C}$, and a non-negative integer $L \in \mathbb{N}$, our goal is to compute approximations $\tilde{y}_j$ for $y_j := F(x_j)$ such that $|\tilde{y}_j - y_j| \leq 2^{-L}$ for all $j = 1, \ldots, n$. Furthermore, let $2^\tau$ and $2^\Gamma$, with $\tau, \Gamma \in \mathbb{N}_{\geq 1}$, denote bounds on the absolute values of the coefficients of $F$ and the points $x_j$, respectively.

For the sake of simplicity, assume that $n = 2^k$ is a power of two; otherwise, pad $F$ with zeros. We require that arbitrarily good approximations of the coefficients $F_i$ and the points $x_j$ are provided by an oracle for the cost of reading the approximations. That is, asking for an approximation of $F$ and the points $x_j$ to a precision of $\ell$ bits after the binary point takes $O(n(\tau + \Gamma + \ell))$ bit operations.

**Algorithm 1 (Multipoint evaluation).** We will follow the classical divide-and-conquer method for fast polynomial multipoint evaluation:

1. From the linear factors $g_{0,j}(x) := x - x_j$, we recursively compute the subproduct tree

$$g_{i,j}(x) := (x - x_{(j-1)2^i+1}) \cdots (x - x_{j2^i}) = g_{i-1,2j-1}(x) \cdot g_{i-1,2j}(x) \tag{1}$$

for $i$ from 1 to $k - 1$ and $j$ from 1 to $n/2^i = 2^{k-i}$, that is, going up from the leaves. Notice that $\deg g_{i,j} = 2^i$.

---

[2] Pan and Tsigaridas [21] claim a comparable complexity bound for their root refinement method when using fast approximate multipoint evaluation. However, they do not provide a rigorous argument to bound the precision demand of their approach.

2. Starting with $r_{k,1}(x) := F(x)$, we recursively compute the remainder tree

$$r_{i,j}(x) := F(x) \bmod g_{i,j}(x) = r_{i+1,\lceil j/2 \rceil}(x) \bmod g_{i,j}(x)$$

for $i$ from $k-1$ to $0$ and $j$ from $1$ to $n/2^i = 2^{k-i}$, that is, going down from the root. Notice that $\deg r_{i,j} < 2^i$.

3. Observe that the value at point $x_j$ is exactly the remainder

$$r_{0,j} = F(x) \bmod g_{0,j}(x) = F(x) \bmod (x - x_j) = F(x_j) \in \mathbb{C}.$$

It is well known that this scheme requires a total number of $O(\mu(n) \log n)$ arithmetic operations in $\mathbb{C}$ (e.g., see [3, Chapter 1, Section 4] or [7, Corollary 10.8]), where $\mu(n)$ denotes the arithmetic complexity of multiplying two polynomials of degree $n$ or, equivalently, the bit complexity of multiplying two $n$-bit integers. Hence, using an asymptotically fast multiplication method with soft-linear bit complexity such as the algorithms by De et al. [4] or Fürer [5] yields a soft-linear arithmetic complexity for polynomial multipoint evaluation. However, we are mainly interested in the bit complexity of the above algorithm if the multiplications and divisions are carried out with approximate but certified arithmetic such that an output precision of $L$ bits after the binary point can be guaranteed. Fast polynomial division is widely considered to be numerically instable which explains why a result on the bit complexity of approximate polynomial evaluation is still missing. We will close this gap by using a method from Schönhage for numerical polynomial division based on a direct application of discrete Fourier transforms to minimize the number of numerically unstable operations; see Section 2.2.

## 2.1. Fast Approximate Polynomial Multiplication

**Definition 2 (Polynomial approximation).** Let $\|\cdot\|$ be a norm on the set of complex polynomials considered as a vector space over $\mathbb{C}$. For a polynomial $f = \sum_{i=0}^n a_i x^i \in \mathbb{C}[x]$ and an integer $\ell$, a polynomial $\tilde{f} \in \mathbb{C}[x]$ is called an *(absolute) $\ell$-bit approximation of $f$ w.r.t.* $\|\cdot\|$ if $\|\tilde{f} - f\| \leq 2^{-\ell}$. Alternatively, if $\tilde{f} = f + \Delta f$, this is equivalent to $\|\Delta f\| \leq 2^{-\ell}$.

When not mentioned explicitly, we assume the norm to be the 1- or sum-norm $\|\cdot\|_1$ with $\|f\|_1 = \sum_{i=0}^n |a_i|$.

The definition of an (absolute) polynomial approximation does not take into account the degree. Typically, degree loss arises when approximating a polynomial with very small leading coefficients which may be truncated to zero. However, the definition also allows for a higher (but finite) degree of the approximation.

We further remark that any $\ell$-bit $\|\cdot\|_1$-approximation of a polynomial implies an $\ell$-bit approximation of each coefficient or, in other words, an $\ell$-bit approximation w.r.t. the $\infty$- or maximum-norm $\|f\|_\infty = \max_i |a_i|$. Conversely, any coefficient-wise approximation $\tilde{f}$ on $f$ to $\ell + \log(\tilde{n} + 1)$ bits, with $\tilde{n} = \deg \tilde{f}$, constitutes an $\ell$-bit $\|\cdot\|_1$-approximation of $f$.

The reason why we favor the sum-norm is its sub-multiplicativity property, that is, for $f, g \in \mathbb{C}[x]$, we have

$$\|f \cdot g\|_1 \leq \|f\|_1 \cdot \|g\|_1. \tag{2}$$

In practice, we also require the precision of the coefficients to be not too high in order to avoid costly arithmetic with superfluous accuracy. Hence, we assume that the coefficients are represented by dyadic values with less than $\ell + \log(n + 1) + c$ bits after the binary point for some small constant $c$.

**Theorem 3 (Numerical multiplication of polynomials).** *Let $f \in \mathbb{C}[x]$ and $g \in \mathbb{C}[x]$ be polynomials of degree less than or equal to $n$ and with coefficients of modulus less than $2^b$ for some integer $b \geq 1$. Then, computing an $\ell$-bit $\|\cdot\|_1$-approximation $\tilde{h}$ for the product $h := f \cdot g$ is possible in*

$$O(\mu(n(\ell + b + 2\log n))) \quad \text{or} \quad \tilde{O}(n(\ell + b))$$

*bit operations and with a precision demand of at most $\ell + O(b + \log n)$ bits on each of the coefficients of $f$ and $g$.*

*Proof.* See, e.g., [7, Corollary 8.27], and note that an overall precision of $\ell + O(b)$ is sufficient to compute an $\tilde{h}$ as desired. For a mostly self-contained description, see Appendix A.1 or [18, Theorem 2.2]. □

## 2.2. Fast Approximate Polynomial Division

**Definition 4 (Numerical division of polynomials).** Given a dividend $f \in \mathbb{C}[x]$, a divisor $g \in \mathbb{C}[x]$, and an integer $\ell \geq 1$, the task of *numerical division of polynomials* is to compute polynomials $\tilde{Q} \in \mathbb{C}[x]$ and $\tilde{R} \in \mathbb{C}[x]$ satisfying

$$\|f - (\tilde{Q} \cdot g + \tilde{R})\|_1 \leq 2^{-\ell}$$

with $\deg \tilde{Q} \leq \deg f - \deg g$ and $\deg \tilde{R} < \deg g$.

**Theorem 5.** *(Schönhage [18, Theorem 4.1]) Let $f \in \mathbb{C}[x]$ be a polynomial of degree $\leq 2n$ and with norm $\|f\|_1 \leq 1$, and let $g \in \mathbb{C}[x]$ be a polynomial of degree $n$ with norm $1 \leq \|g\|_1 \leq 2$. Suppose that a bound $2^\rho$, with $\rho \in \mathbb{N}_{\geq 1}$, on the modulus of all roots of $g$ is given. Then, numerical division of $f$ by $g$ up to an error of $2^{-\ell}$ needs a number of bit operations bounded by*

$$O(\mu(n(\ell + n\rho))) = \tilde{O}(n(\ell + n\rho)).$$

In his presentation of the division algorithm, Schönhage carefully analyses the required precision for the needed operations in his algorithm as $2 \cdot \ell + 5n\rho + O(n)$ bits; see [18, (4.14) and (4.15)]. Hence, one might conclude that this bound also expresses the precision demand on the input polynomials $f$ and $g$. However, the factor **2** in the above bound is only needed for the precision with which the internal computations have to be performed, whereas it is not necessary for the precision demand of the input polynomials $f$ and $g$. In particular, for $\ell \ggg n\rho$, input and output accuracy are asymptotically identical, independently from the algorithm used to carry out the numerical division. For a proof of the above claim, we need an additional result from Schönhage which provides a worst-case perturbation bound for polynomial zeros under perturbation of its coefficients.

**Theorem 6.** *(Schönhage [19, Theorem 2.7]) Let $f \in \mathbb{C}[x]$ be a polynomial of degree $n$ with zeros $x_1, \ldots, x_n$, not necessarily distinct, and let $\hat{f}$ be a $\log(\eta \|f\|_1)$-approximation of $f$ for $\eta \leq 2^{-7n}$. Then, the zeros $\hat{x}_1, \ldots, \hat{x}_n$ of $\hat{f}$ can be numbered such that $|\hat{x}_j - x_j| < 9 \sqrt[n]{\eta}$ for $|x_j| \leq 1$ and $|1/\hat{x}_j - 1/x_j| < 9 \sqrt[n]{\eta}$ for $|x_j| \geq 1$.*[3]

We can now give a stronger version of Theorem 5 which comprises the claimed bound on the needed input precision. In addition, we show that, within a comparable time bound as given in Theorem 5, we can guarantee that the computed polynomials $\tilde{Q}$ and $\tilde{R}$ are $\ell$-bit approximations of their exact counterparts.

**Theorem 7.** *Let $f$, $g$ and $\rho$ as in Theorem 5, and $Q := f \operatorname{div} g$ and $R := f \bmod g$ be the exact quotient and remainder in the polynomial division of $f$ by $g$.*

*Then, the cost for computing $\ell$-bit approximations $\tilde{Q}$ and $\tilde{R}$ of $Q$ and $R$ satisfying $\|f - (\tilde{Q} \cdot g + \tilde{R})\|_1 \leq 2^{-\ell}$ is bounded by $\tilde{O}(n(\ell + n\rho))$ bit operations. For this computation, we need $(\ell + 32n\rho)$-bit approximations of the polynomials $f$ and $g$.*

---

[3] Schönhage points out that the theorem also holds for zeros at infinity, that is, in the case where $\deg f \neq \deg \hat{f}$. However, in our applications, the degrees will always be the same.

*Proof.* Let $\hat{f} = f + \Delta f$ and $\hat{g} = g + \Delta g$ be arbitrary $\ell_f$- and $\ell_g$-bit approximations for $f = \sum_{i=0}^{2n} f_i x^i$ and $g = \sum_{i=0}^n g_i x^i$, where $\deg \hat{f} \le 2n$, $\deg \hat{g} \le \deg g$, and $\ell_f$ and $\ell_g$ are integers to be specified later.

First, we note that $\deg g$ and $\deg \hat{g}$ actually coincide for any $\ell_g \ge n(\rho+2)+1$. Namely, there exists at least one coefficient $g_i$ of $g$ with $|g_i| \ge 1/(n+1) \ge 2^{-n}$ since $\|g\|_1 \ge 1$, and thus $|g_n| \ge |g_i| \cdot 2^{-n-n\rho} \ge 2^{-n(\rho+2)}$, where the second to last inequality follows from the fact that $|g_i| \le |g_n| \cdot 2^n \prod_{z:g(z)=0} \max(1,|z|) \le |g_n| \cdot 2^{n+n\rho}$. Hence, in particular, we have

$$|g_n|, |\hat{g}_n| \ge 2^{-n(\rho+2)-1} \ge 2^{-4n\rho} \text{ for all } \ell_g \ge n(\rho+2)+1.$$

Next, we derive a necessary condition on the precision $\ell_g$ such that $2^{2\rho}$ is a root bound for $\hat{g}$. Suppose that $\ell_g \ge \max(n(\rho+2)+1, 7n+1)$, then we may apply Theorem 6 to the polynomials $g$ and $\hat{g}$ which have the same degrees as shown above. For $x$ and $\hat{x}$ an arbitrary corresponding pair of roots of the polynomials $g$ and $\hat{g}$, we distinguish two cases:

1. If $|x| \le 1$, it immediately follows $|\hat{x}| < |x| + 9\sqrt[n]{2^{-\ell_g}}$ and, hence, $|\hat{x}| < 2^\rho + 1 < 2^{2\rho}$.

2. For $x$ outside the unit circle, we have $|x|/|\hat{x}| > 1 - 9\sqrt[n]{2^{-\ell_g}}|x|$. Thus, we aim for $9\sqrt[n]{2^{-\ell_g}} 2^\rho \le 1/2$ which is fulfilled if $\ell_g \ge n(\rho+5) > n \log 18 + n\rho$.

In what follows, we assume that $\ell_g \ge \max(7n+1, n(\rho+5))$. This ensures that the degrees of $g$ and $\hat{g}$ coincide and $2^{\hat{\rho}}$, with $\hat{\rho} := 2\rho$, constitutes an upper bound on the absolute value of all roots of $g$ as well as for all roots of $\hat{g}$.

Suppose that $f = Q \cdot g + R$ and $\hat{f} = \hat{Q} \cdot \hat{g} + \hat{R}$ are the *exact* representations of $f$ and $\hat{f}$ after division with remainder, then we aim to show that the pair $(\hat{Q}, \hat{R})$ is actually a good approximation of $(Q, R)$ (i.e., $\approx \min(\ell_f, \ell_g)$-bit approximations) if $\ell_f$ and $\ell_g$ are both large enough. Write $\Delta Q := \hat{Q} - Q$ and $\Delta R := \hat{R} - R$. The coefficients $Q_k$ of $Q$ appear as leading coefficients in the Laurent series of the function

$$\frac{f(x)/x^n}{g(x)} = \frac{f_{2n} + f_{2n-1}/x + f_{2n-2}/x^2 + \cdots}{g_n + g_{n-1}/x + g_{n-2}/x^2 + \cdots} = Q_n + \frac{Q_{n-1}}{x} + \frac{Q_{n-2}}{x^2} + \cdots$$

and can be represented, using Cauchy's integral formula, as

$$Q_k = \frac{1}{2\pi i} \int_{|x|=\varrho} \frac{f(x)/x^n}{g(x)} x^{k-1} \, dx \qquad (3)$$

for any $\varrho > 2^\rho$; see [18, (4.7)–(4.9)]. Using the corresponding representation of the coefficients $\hat{Q}_k$ of $\hat{Q}$, we can estimate (here, for any $\varrho > 2^{\hat{\rho}}$)

$$|\hat{Q}_k - Q_k| = \frac{1}{2\pi} \left| \int_{|x|=\varrho} \frac{\Delta f(x) \cdot g(x) - f(x) \cdot \Delta g(x)}{g(x) \hat{g}(x)} x^{k-n-1} \, dx \right|. \qquad (4)$$

Throughout the following considerations, we fix $\varrho := 2^{\hat{\rho}} + 1 = 2^{2\rho} + 1 < 2^{3\rho}$. The absolute value of the numerator of the integrand is bounded by

$$(|\Delta f(x) \cdot g(x)| + |f(x) \cdot \Delta g(x)|) \cdot |x|^{k-n-1}$$
$$\le (\|\Delta f\|_1 \varrho^{2n} \cdot \|g\|_1 \varrho^n + \|f\|_1 \varrho^{2n} \cdot \|\Delta g\|_1 \varrho^n) \cdot \varrho^{k-n-1}$$
$$\le (2^{-\ell_f+1} + 2^{-\ell_g+1}) \varrho^{2n+k-1} \le (2^{-\ell_f+1} + 2^{-\ell_g+1}) \varrho^{2n+k}$$

and, for the denominator, we have

$$|g(x)\,\hat{g}(x)| \geq |g_n|(\varrho - 2^\rho)^n \cdot |\hat{g}_n|(\varrho - 2^{\hat{\rho}})^n \geq |g_n| \cdot |\hat{g}_n| \geq 2^{-8n\rho}.$$

Now, using the latter two estimates in (4) yields

$$|\Delta Q_k| = |\hat{Q}_k - Q_k| \leq (2^{1-\ell_g} + 2^{-\ell_g}) \cdot 2^{8n\rho} \cdot \varrho^{2n+k}.$$

Summing over all $k = 0, \ldots, n$ gives

$$\begin{aligned}
\|\Delta Q\|_1 = \|\hat{Q} - Q\|_1 &\leq (2^{-\ell_f+1} + 2^{-\ell_g+1}) \cdot 2^{8n\rho} \cdot \varrho^{2n} \frac{\varrho^{n+1} - 1}{\varrho - 1} \\
&\leq (2^{-\ell_f+1} + 2^{-\ell_g+1}) \cdot 2^{8n\rho} \cdot \varrho^{3n+1} \\
&< (2^{-\ell_f+1} + 2^{-\ell_g+1}) \cdot 2^{20n\rho} < 2^{-\min(\ell_f,\ell_g)+20n\rho+2},
\end{aligned}$$

where we used that $\varrho = 2^{2\rho} + 1 < 2^{2\rho+1}$ and thus $\varrho^{3n+1} < 2^{12n\rho}$. Hence, for

$$\ell_f, \ell_g \geq \ell + 20n\rho + 4, \tag{5}$$

the polynomial $\hat{Q}$ is an $(\ell + 2)$-bit approximation of $Q$. An analogous computation as above based on the formula (3) further shows that

$$\|Q\|_1 \leq 2^{4n\rho} \cdot \varrho^{2n+1} \leq 2^{13n\rho}. \tag{6}$$

Hence, under the above constraints from (5) for $\ell_f$ and $\ell_g$, we conclude that

$$\begin{aligned}
\|\hat{R} - R\|_1 &= \|(\hat{f} - \hat{Q}\,\hat{g}) - (f - Q\,g)\|_1 \\
&\leq \|\Delta f\|_1 + \|Q\|_1 \|\Delta g\|_1 + \|\Delta Q\|_1 \|g\|_1 + \|\Delta Q\|_1 \|\Delta g\|_1 \\
&\leq 2^{-\ell_f} + 2^{13n\rho} \cdot 2^{-\ell_g} + 2^{-\ell-2} \cdot 2 + 2^{-\ell-2} \cdot 2^{-\ell_g} \\
&< 2^{-\ell-3} + 2^{-\ell-3} + 2^{-\ell-1} + 2^{-\ell-3} \leq 2^{-\ell},
\end{aligned}$$

thus $\hat{R}$ constitutes an $\ell$-bit approximation of $R$.

We are now in the position to put the pieces together and prove the main statements of the theorem. For $\tilde{\ell} := \ell + 32n\rho > (\ell + 3) + 20n\rho + 4$,[4] we first choose $\tilde{\ell}$-bit approximations $\tilde{f}$ and $\tilde{g}$ of $f$ and $g$, respectively, such that $\|\tilde{f}\|_1 \leq 1$ and $1 \leq \|\tilde{g}\|_1 \leq 2$.[5] We can now apply Theorem 5 to compute polynomials $\tilde{Q}$ and $\tilde{R}$ such that $\|\tilde{f} - (\tilde{Q} \cdot \tilde{g} + \tilde{R})\|_1 \leq 2^{-\tilde{\ell}}$. For this step, we need $\tilde{O}(n(\ell + n\rho))$ bit operations. We define $\hat{f} := \tilde{Q} \cdot \hat{g} + \tilde{R}$ and $\hat{g} := \tilde{g}$, where the latter two polynomials are $\tilde{\ell}$-bit approximations of $f$ and $g$, respectively. Thus, the above consideration shows that $\tilde{Q}$ and $\tilde{R}$ are $(\ell + 3)$-bit approximations of the exact solutions $Q$ and $R$, respectively. It follows that

$$\|f - (\tilde{Q}g + \tilde{R})\|_1 \leq \|(Q - \tilde{Q}) \cdot g\|_1 + \|R - \tilde{R}\|_1 \leq \|Q - \tilde{Q}\|_1 \cdot \|g\|_1 + 2^{-\ell-3} \leq 2^{-\ell}$$

holds, completing the proof. □

---

[4] In fact, it suffices to choose $\tilde{\ell} := \ell + 27n\rho$, however, we aimed for "nice numbers."
[5] Notice that this can always be achieved since we can always choose approximations of $f$ and $g$ which decrease or increase the corresponding 1-norms by less than $2^{-\tilde{\ell}} < 1/2$.

Notice that the above result shows that the precision demand for the input polynomials is of the same size as the desired output precision plus a term which only depends on fixed parameters, that is, $n$ and $\rho$. This will turn out to be crucial when considering numerical division within the multipoint evaluation algorithm. Namely, since we have to perform $\log n$ successive divisions, a precision demand of $2 \cdot \ell$ (as needed for the internal computations in Schönhage's algorithm) for the input in each iteration would eventually propagate to a precision demand of $n\ell$, which is undesirable. However, from the above theorem, we conclude that, for an output precision of $\ell$, an input precision of $\ell + O((\log n) \cdot n\rho)$ is sufficient because, in each of the $\log n$ successive divisions, we loose a precision of $O(n\rho)$. In order to give more precise results (and rigorous arguments), we first have to make Theorem 7 applicable to polynomials with higher norm as they appear in the multipoint evaluation scheme. With this task in mind, we concentrate on the case of *monic* divisors $g$.

**Corollary 8.** *Let $f \in \mathbb{C}[x]$ be a complex polynomial of degree $\leq 2n$ with $\|f\|_1 \leq 2^b$ for some integer $b \geq 1$, and let $g$ be a* monic *polynomial of degree $n$ with a given root bound $2^\rho$, where $\rho \in \mathbb{N}_{\geq 1}$. Let $Q := f \operatorname{div} g$ and $R := f \bmod g$ denote the exact quotient and remainder in the polynomial division of $f$ by $g$.*

*Then, the cost for computing $\ell$-bit approximations $\tilde{Q}$ and $\tilde{R}$ of $Q$ and $R$, respectively, with $\|f - (\tilde{Q} \cdot g + \tilde{R})\|_1 \leq 2^{-\ell}$ is bounded by $\tilde{O}(n(\ell + b + n\rho))$ bit operations. For this computation, we need $(\ell + b + n(2\rho + 2\lceil \log 2n \rceil + 32))$-bit approximations of the polynomials $f$ and $g$. The approximate remainder $\tilde{R}$ fulfills*

$$\|\tilde{R}\|_1 \leq 2^{16n + 2n\rho + 2n\log\lceil 2n\rceil + b} = 2^{b + 2n\rho + O(n\log n)}. \tag{7}$$

*Proof.* Let $s := \rho + \lceil \log 2n \rceil$ and $\ell^* := \ell + b + ns$. We define

$$f^*(x) := 2^{-b-2ns} f(2^s x) \quad \text{and} \quad g^*(x) := 2^{-ns} g(2^s x).$$

It follows that $\|f^*\|_1 \leq 1$ and $1 \leq \|g^*\|_1$ since $g^*$ is again monic. In addition, the scaling of $g$ by $2^s$ yields that all roots of $g^*$ have absolute values less than or equal to $1/(2n)$. Thus, the $i$-th coefficient of $g^*$ is bounded by $\binom{n}{i}\frac{1}{(2n)^i}$ which shows that $\|g^*\|_1 \leq (1 + \frac{1}{2n})^n < e^{1/2} < 2$. We now apply Theorem 7 to the polynomials $f^*$ and $g^*$, and to some desired output precision $\ell^*$ which will be specified later: Suppose that $Q^* := f^* \operatorname{div} g^*$ and $R^* := f^* \bmod g^*$, it takes $\tilde{O}(n(\ell^* + n))$ bit operations to compute $\ell^*$-bit operations $\tilde{Q}^*$ and $\tilde{R}^*$ of $Q^*$ and $R^*$, respectively, such that $\|f^* - (\tilde{Q}^* g + \tilde{R}^*)\|_1 < 2^{-\ell^*}$. For this, we need $(\ell^* + 32n)$-bit operations of the polynomials $f^*$ and $g^*$. Notice that we used the fact that $2^1$ constitutes a root bound for $g^*$. We further remark that, in order to compute the approximations for $f^*$ and $g^*$, it suffices to consider $(\ell^* + 32n)$-bit approximations of the polynomials $f$ and $g$. In order to recover approximations for the polynomials $Q$ and $R$, we consider an inverse scaling, that is,

$$\tilde{Q}(x) := 2^{b+ns} \tilde{Q}^*(2^{-s} x) \quad \text{and} \quad \tilde{R}(x) := 2^{2ns+b} \tilde{R}^*(2^{-s} x).$$

Since $f(x) = Q(x) \cdot g(x) + R(x)$, we have

$$\underbrace{2^{-2ns-b} f(2^s x)}_{f^*} = \underbrace{2^{-b-ns} Q(2^s x)}_{Q^*} \cdot \underbrace{2^{-ns} g(2^s x)}_{g^*} + \underbrace{2^{-2ns-b} R(2^s x)}_{R^*},$$

and, thus, for any $\ell^* \geq b + 2ns$, the polynomials $\tilde{Q}(x)$ and $\tilde{R}(x)$ are $(\ell^* - b - 2ns)$-approximations of $Q$ and $R$, respectively. In addition, $\|f - (\tilde{Q}g + \tilde{R})\|_1 \leq 2^{-\ell^* + b + 2ns}$. Hence, for $\ell^* := \ell + b + 2ns$, the bound on the bit complexity of the numerical division as well as the bound on the precision demand follows.

For the estimate on $\|\tilde{R}\|_1$, we recall that (6) yields $\|Q^*\|_1 \leq 2^{13n}$ which implies that $\|R^*\|_1 \leq \|f^*\|_1 + \|Q^*\|_1 \cdot \|g^*\|_1 \leq 1 + 2 \cdot 2^{13n} < 2^{16n}$. Thus, we have

$$\|R\|_1 \leq 2^{2ns+b} \|R^*\|_1 < 2^{16n + 2ns + b} < 2^{16n + 2n\rho + 2n\log\lceil 2n\rceil + b} = 2^{b + 2n\rho + O(n\log n)},$$

and the same bound also applies to $\tilde{R}$ since it is an $\ell$-bit approximation of $R$. $\square$

## 2.3. Fast Approximate Multipoint Evaluation: Complexity Analysis

We can now apply the results of the previous two sections to the recursive divide-and-conquer multipoint evaluation scheme as described on page 3. Using approximate multiplications and divisions, our goal is to compute approximations $\tilde{r}_{0,j}$ of the final remainders $r_{0,j} = F \bmod (x - x_j) = F(x_j)$ such that $|\tilde{r}_{0,j} - F(x_j)| \leq 2^{-L}$ for all $j = 1, \ldots, n$. In other words, we aim to compute $L$-bit approximations of the remainders $r_{0,j}$. For this purpose, we will do a bottom-up backwards analysis of the required precisions for dividend and divisor in every layer of the remainder tree which will yield the according requirements on the accuracy of the subproduct tree.

**Theorem 9.** *Let $F \in \mathbb{C}[x]$ be a polynomial of degree $n$ with $\|F\|_1 \leq 2^\tau$, with $\tau \geq 1$, and let $x_1, \ldots, x_n \in \mathbb{C}$ be complex points with absolute values bounded by $2^\Gamma$, where $\Gamma \geq 1$. Then, approximate multipoint evaluation up to a precision of $2^{-L}$ for some integer $L \geq 0$, that is, computing $\tilde{y}_j$ such that $|\tilde{y}_j - F(x_j)| \leq 2^{-L}$ for all $j$, is possible with*

$$\tilde{O}(n(L + \tau + n\Gamma)).$$

*bit operations. Moreover, the precision demand on $F$ and the points $x_j$ is bounded by $L + O(\tau + n\Gamma + n \log n)$ bits.*

*Proof.* Define $g_{i,j}$ and $r_{i,j}$ as in Algorithm 1. We analyse a run of the algorithm using approximate multiplication and division, with a precision of $\ell_i^{\mathrm{div}}$ for the approximate divisors $\tilde{g}_{i,*}$ and remainders $\tilde{r}_{i,*}$ in the $i$-th layer of the subproduct and the remainder tree. We recall that $\deg \tilde{g}_{i,*} = \deg g_{i,*} = 2^i$.

According to Corollary 8, for the recursive divisions to yield an output precision $\ell_i \geq 0$, it suffices to have approximations $\tilde{r}_{i+1,*}$ and $\tilde{g}_{i,*}$ of the exact polynomials $f := r_{i+1,*}$ and $g := g_{i,*}$ to a precision of

$$\ell_{i+1}^{\mathrm{div}} := \ell_i^{\mathrm{div}} + \log \|r_{i+1,*}\|_1 + 2^{i+1}\Gamma + O(i \cdot 2^i) \tag{8}$$

bits, since the roots of each $g_{i,*}$ are contained in the set $\{x_1, \ldots, x_n\}$ and, thus, their absolute values are also bounded by $2^\Gamma$. In addition, it holds that $\|r_{\log n, 0}\|_1 = \|F\|_1 \leq 2^\tau$. In order to bound the absolute values of the remainders $r_{i,*}$ for $i < \log n$, we can use our remainder bound from (7) in an iterative manner to show that

$$\log \|r_{i,*}\|_1 = \log \|r_{i+1,*}\|_1 + 2^{i+1}\Gamma + O(i \cdot 2^i) = \tau + 2n\Gamma + O(n \log n). \tag{9}$$

Combining (8) and (9) then yields $\ell^{\mathrm{div}} := \max_{i > 0} \ell_i^{\mathrm{div}} = \ell_0^{\mathrm{div}} + \tau + 2n\Gamma + O(n \log n)$. Hence, choosing $\ell_0^{\mathrm{div}} := L$, we eventually achieve evaluation up to an error of $2^{-L}$ if all numerical divisions are carried out with precision $\ell^{\mathrm{div}}$. The bit complexity to carry out a single numerical division at the $i$-th layer of the tree is then bounded by $\tilde{O}(2^i(\ell^{\mathrm{div}} + \tau + 2^i\Gamma)) = \tilde{O}(2^i(L + n\Gamma + \tau))$. Since there are $n/2^i$ divisions, the total cost at the $i$-th layer is bounded by $\tilde{O}(n(L + n\Gamma + \tau))$. The depth of the tree equals $\log n$, and thus the overall bit complexity is $\tilde{O}(n(L + n\Gamma + \tau))$.

It remains to bound the precision demand and the cost for computing $(L + \tau + 2n\Gamma + O(n \log n))$-bit approximations of the polynomials $g_{i,*}$. According to Theorem 3, in order to compute the polynomials $g_{i,*}$ to a precision of $\ell_i^{\mathrm{mul}}$, we have to consider $\ell_{i-1}^{\mathrm{mul}}$-bit approximations of $g_{i-1,*}$, where $\ell_i^{\mathrm{mul}} := \ell_{i-1}^{\mathrm{mul}} + 2\log \|g_{i-1,*}\|_1 + O(i) = \ell_{i-1}^{\mathrm{mul}} + 2i\Gamma + O(i) = \ell_0^{\mathrm{mul}} + O(\log n \cdot \Gamma)$. Hence, it suffices to run all multiplications in the product tree with a precision of $\ell^{\mathrm{mul}} = L + \tau + O(n\Gamma + n \log n)$. The bit complexity for all multiplications is bounded by $\tilde{O}(n\ell^{\mathrm{mul}}) = \tilde{O}(n(L + \tau + n\Gamma))$, and the precision demand for the points $x_i$ is bounded by $\ell^{\mathrm{mul}} + O(\Gamma + \log n) = L + \tau + O(n\Gamma + n \log n)$. □

## 3. Applications

### 3.1. Polynomial Interpolation

Fast polynomial interpolation can be considered as a direct application of polynomial multipoint evaluation. Given $n$ (w.l.o.g. we again assume that $n = 2^k$ is a power of two) pairwise distinct interpolation points $x_1, \ldots, x_n \in \mathbb{C}$ and corresponding values $v_1, \ldots, v_n$, we aim to compute the unique polynomial $F \in \mathbb{C}[x]$ of degree less than $n$ such that $F(x_i) = v_i$ for all $i = 1, \ldots, n$. Using Lagrange interpolation, we have

$$F(x) = \sum_{i=1}^{n} v_i \cdot \prod_{j=1;\, j \neq i}^{n} \frac{x - x_j}{x_i - x_j} = \sum_{i=1}^{n} v_i \lambda_i^{-1} \cdot \prod_{j=1;\, j \neq i}^{n} (x - x_j) = \sum_{i=1}^{n} \mu_i \cdot \prod_{j=1;\, j \neq i}^{n} (x - x_j),$$

where $\lambda_i := \prod_{j=1;\, j \neq i}^{n}(x_i - x_j)$ and $\mu_i := v_i \cdot \lambda_i^{-1}$. Now, in order to compute $F(x)$, we proceed in two steps: In the first step, we compute the values $\lambda_i$. Let $g(x) := \prod_{j=1}^{n}(x - x_j)$ (notice that $g(x)$ coincides with the polynomial $g_{k,1}(x)$ from (1)), then $\lambda_i = g'(x_i)$, and thus the values $\lambda_i$ can be obtained by a fast multipoint evaluation of the derivative $g'(x)$ of the polynomial $g(x)$ at the points $x_i$. We can compute $g$ and $g'$ with $\tilde{O}(n)$ arithmetic operations in $\mathbb{C}$, and, using fast multipoint evaluation, the same bound also applies to the number of arithmetic operations to compute all values $\lambda_i$. Hence, computing the values $\mu_i$ takes $\tilde{O}(n)$ arithmetic operations in $\mathbb{C}$. Now, in order to compute $F_{k,1}(x) := F(x) = \sum_{i=1}^{n} \mu_i \cdot \prod_{j=1;\, j \neq i}^{n}(x - x_j)$, we write

$$F_{k,1}(x) = g_{k-1,1}(x) \cdot \underbrace{\sum_{i=1}^{n/2} \mu_i \cdot \prod_{j=1;\, j \neq i}^{n/2}(x - x_j)}_{=: F_{k-1,1}(x)} + g_{k-1,2}(x) \cdot \underbrace{\sum_{i=n/2+1}^{n} \mu_i \cdot \prod_{j=n/2+1;\, j \neq i}^{n}(x - x_j)}_{=: F_{k-1,2}(x)}. \tag{10}$$

Following a divide-and-conquer approach, we can recursively compute $F(x)$ from the values $\mu_i$ and the polynomials $g_{i,j}$ as defined in (1). It is then straight forward to show that $\tilde{O}(n)$ arithmetic operations in $\mathbb{C}$ are sufficient to carry out the necessary computations.

In contrast to the exact computation of $F(x)$ as outlined above, we now focus on computing an $L$-bit approximation $\tilde{F}$ of $F$. We assume that arbitrarily good approximations of the points $x_i$ and the corresponding values $v_i$ are provided. We introduce the following definitions:

$$\begin{aligned}\Gamma &:= \max_{i=1}^{n} \log \max(2, |x_i|) \geq 1, \quad V := \max_{i=1}^{n} \log \max(2, |v_i|) \geq 1, \quad \text{and} \\ \Lambda &:= \max_{i=1}^{n} \log \max(1, |\lambda_i|^{-1}) = \max_{i=1}^{n} \log \max(1, \prod_{j=1;\, j \neq i}^{n} |x_i - x_j|^{-1}).\end{aligned} \tag{11}$$

In Section 2.1, we have already shown that computing $\ell$-bit approximations of all $g_{i,j}$ needs $\tilde{O}(n^2\Gamma + n\ell)$ bit operations. Furthermore, we need approximations of the points $x_i$ to $\tilde{O}(n\Gamma + \ell)$ bits after the binary point. Applying Theorem 9 to the derivative $g'(x) := g'_{1,k}(x)$ of $g_{1,k}(x)$ and the points $x_i$ then shows that computing $\ell$-bit approximations $\tilde{\lambda}_i$ of the values $\lambda_i$ uses $\tilde{O}(n^2\Gamma + n\ell)$ bit operations since the modulus of all $x_i$ is bounded by $2^\Gamma$ and the coefficients of $g'$ have absolute value $2^{O(n\Gamma)}$. The precision demand on $g'$ and the points $x_i$ is bounded by $\tilde{O}(n\Gamma + \ell)$ bits after the binary point. Now, in order to compute an $\ell$-bit approximation $\tilde{\mu}_i$ of $\mu_i = v_i \lambda_i^{-1}$, we have to approximate $v_i$ and $\lambda_i$ to $O(\ell + \Lambda + V)$ bits. Hence, computing such approximations $\tilde{\mu}_i$ for all $i$ needs $\tilde{O}(n(\ell + n\Gamma + \Lambda + V))$ bit operations, and the precision demand for the points $x_i$ and the values $v_i$ is bounded by $\tilde{O}(\ell + n\Gamma + \Lambda + V)$ bits. For computing $\tilde{F}$, we now apply the recursion from (10). Starting with $\ell$-bit approximations of $\mu_i$ and $g_{i,j}$, the so-obtained polynomial $\tilde{F}$ differs from $F$ by at most $\ell - O(n\Gamma + \Lambda + V)$-bits after the binary point since the coefficients of all occurring polynomials in the intermediate computations have modulus bounded by $2^{O(n\Gamma + \Lambda + V)}$. Hence, we conclude the following theorem:

**Theorem 10.** *Let $x_1, \ldots, x_n \in \mathbb{C}$ be arbitrary, but distinct, given interpolation points and $v_1, \ldots, v_n \in \mathbb{C}$ be arbitrary corresponding interpolation values. Furthermore, let $F \in \mathbb{C}[x]$ be the unique polynomial of degree less than $n$ such that $F(x_i) = v_i$ for all $i$. Then, for any given integer $L$, we can compute an $L$-bit approximation $\tilde{F}$ of $F$ with*

$$\tilde{O}(n(n\Gamma + V + \Lambda + L)) \tag{12}$$

*bit operations, where $\Gamma$, $V$, and $\Lambda$ are defined as in* (11). *The points $x_i$ and the values $v_i$ have to approximated to $\tilde{O}(n\Gamma + V + \Lambda + L)$ bits after the binary point.*

*Remark* 11. In the special case, where $x_i = e^{i \cdot \frac{2\pi i}{n}}$ are the $n$-th roots of unity, we have $\Gamma = 1$ and $\Lambda = \log n$ because $\prod_{j=1; j \neq i}^{n} |x_i - x_j| = \left|\frac{d(x^n - 1)}{dx}(x_i)\right| = |n \cdot x_i^{n-1}| = n$. The bound in (12) then simplifies to $\tilde{O}(n(n + L + V))$ which is comparable to the complexity bound that one gets from considering an inverse FFT to the vector $(v_1, \ldots, v_n)$ using approximate arithmetic [20, Theorem 8.3], regardless of the fact that the latter approach is certainly more reasonable.

### 3.2. Asymptotically Fast Approximate Taylor Shifts

Our second application concerns the problem of computing the Taylor shift of a polynomial $F \in \mathbb{C}[x]$ by a given $m \in \mathbb{C}$. More precisely, given oracles for arbitrarily good approximations of $F$ and $m$ and a positive integer $L$, we aim to compute an $L$-bit approximation of $F_m(x) := F(m + x)$. Computing the shifted polynomial $F_m$ is crucial in many subdivision algorithms to compute the roots of a polynomial. Asymptotically fast methods have already been studied in [20] and [6], where the computation of the coefficients of $F_m$ is reduced to a multiplication of two polynomials. We follow a slightly different approach based on multipoint evaluation, where the problem is reduced to an evaluation-interpolation problem. More specifically, we first evaluate $F$ at the $n$ points $x_i := m + e^{i \cdot \frac{2\pi i}{n}}$, where $n := \deg F + 1$. We then compute $F_m$ as the unique polynomial of degree less than $n$ which takes the values $v_i := p(x_i)$ at the roots of unity $\omega_i := e^{i \cdot \frac{2\pi i}{n}}$. In the preceding sections, we have shown how to carry out the latter two computations with an output precision of $\ell$ bits after the binary point. Theorem 10 and the subsequent remark shows that, in order to compute an $L$-bit approximation of $F_m$, it suffices to run the final interpolation with an input precision of $\tilde{O}(n + L + V)$ bits after the binary point, where $V = \max_{i=1}^{n} \log \max(2, |F(x_i)|) = \log n 2^{\tau} (2 \max(1, |m|))^n = O(n + \tau + n \log \max(1, |m|))$ and $\|F\|_\infty < 2^\tau$. The cost for the interpolation is bounded by $\tilde{O}(n(n + L + V)) = \tilde{O}(n(n + L + \tau + n \log \max(1, |m|)))$. It remains to bound the cost for the evaluation of $F$ at the points $x_i$. Since we need approximations of $F(x_i)$ to $\tilde{O}(L + n + \tau + n \log \max(1, |m|))$ bits after the binary point and $|x_i| < 2 \max(1, |m|)$ for all $i$, Theorem 9 yields that we need $\tilde{O}(n(n + n \log \max(1, |m|) + \tau + L))$ bit operations to run the approximate multipoint evaluation. The polynomial $F$ and the points $x_i$ have to be approximated to $\tilde{O}(n + n \log \max(1, |m|) + \tau + L)$ bits after the binary point. We fix the following result which provides a complexity bound comparable to [20, Theorem 8.4]:

**Theorem 12.** *Let $F \in \mathbb{C}[x]$ be a polynomial of degree less than $n$ with coefficients of modulus less than $2^\tau$, and let $m \in \mathbb{C}$ be an arbitrary complex number. Then, for any $L \in \mathbb{N}_{\geq 1}$, we can compute an $L$-bit approximation $\tilde{F}_m$ of $F_m(x) = F(m + x)$ with $\tilde{O}(n(n + n \log \max(1, |m|) + \tau + L))$ bit operations. For this computation, the coefficients of the polynomial $F$ as well as the point $m$ have to be approximated to $\tilde{O}(n + n \log \max(1, |m|) + \tau + L)$ bits after the binary point.*

### 3.3. Quadratic Interval Refinement for Roots of a Polynomial

As a last application, we very briefly sketch how multipoint evaluation can be used in root finding and refinement. Classical approaches for real root isolation on square-free polynomials $F(x) \in \mathbb{Z}[x]$, such as variants of the Descartes method [17], start from an initial interval $I_0 \subset \mathbb{R}$ known to comprise all roots of $F$ and subdivide $I_0$ to eventually yield disjoint intervals $I_\alpha$, each containing a real root $\alpha$ of $F$. If more precise approximations are required, it conceptually suffices to iterate the same process further until the intervals shrink below the desired threshold.

In [1], Abbott gave a scheme that improves upon the most simple refinement scheme of interval bisection. He proposes a combination of bisection and a secant method to eventually achieve quadratic convergence to the roots, given a set of already isolating intervals. This approach has been adapted in [9, 10] for approximate arithmetic. It exclusively relies on polynomial evaluations on points in the isolating intervals. Our multipoint scheme allows to perform those computations for *all* roots simultaneously with only a polylogarithmic overhead compared to only one evaluation. Hence, we can achieve the following complexity result for real root isolation and refinement. (For a more detailed description, see Appendix A.2 or [10].)

**Theorem 13.** *Let $F \in \mathbb{Z}[x]$ be a square-free polynomial of degree n with integer coefficients bounded by $2^\tau$, and let L be an arbitrary given positive integer. Then, computing isolating intervals for all real roots of F of width $2^{-L}$ or less uses $\tilde{O}(n^3\tau + nL)$ bit operations.*

# A. Appendix

## A.1. Proof of Theorem 3

**Definition 14 (Integer truncation).** For a complex number $z = a + \mathrm{i}\, b \in \mathbb{C}$, a Gaussian integer $\tilde{z} = \tilde{a} + \mathrm{i}\, \tilde{b} \in \mathbb{Z}[\mathrm{i}]$ is called an *integer truncation of* $z$ if $|z - \tilde{z}| \leq 1$. An integer truncation $\tilde{f} \in \mathbb{Z}[\mathrm{i}][x]$ of a polynomial $f \in \mathbb{C}[x]$ is defined coefficient-wise.

In what follows, we ignore the fact that there are several truncations for a complex number and, for the sake of simplicity, pretend that we can compute "the" truncation of $f$ and denote it by $\mathrm{trunc}(f)$. This is reasonable since, for any $\ell \geq 1$, coefficient-wise rounding of any $\ell$-bit $\|\cdot\|_\infty$-approximation of $f$ yields a unique truncation (although not necessarily the same for different $\ell$).

Let us now recall Theorem 3:

**Theorem 3 (Numerical multiplication of polynomials).** *Let $f \in \mathbb{C}[x]$ and $g \in \mathbb{C}[x]$ be polynomials of degree less than or equal to $n$ and with coefficients of modulus less than $2^b$ for some integer $b \geq 1$. Then, computing an $\ell$-bit $\|\cdot\|_1$-approximation $\tilde{h}$ for the product $h := f \cdot g$ is possible in*

$$O(\mu(n(\ell + b + 2\log n))) \quad \text{or} \quad \tilde{O}(n(\ell + b))$$

*bit operations and with a precision demand of at most $\ell + O(b + \log n)$ bits on each of the coefficients of $f$ and $g$.*

*Proof.* Let $s := \ell + b + 2\lceil \log(n+1) \rceil + 2$. Define $F := 2^s f$ and $G := 2^s g$, and notice that $H := FG = 2^{2s} h$. We consider polynomials $\tilde{F} := \mathrm{trunc}(F)$ and $\tilde{G} := \mathrm{trunc}(G) \in \mathbb{Z}[\mathrm{i}][x]$ and write $\Delta F := \tilde{F} - F$ and $\Delta G := \tilde{G} - G$. Since $\|\Delta F\|_1, \|\Delta G\|_1 \leq n+1$,

$$\begin{aligned}
\|\tilde{F}\tilde{G} - FG\|_1 &\leq \|\Delta F \cdot G\|_1 + \|F \cdot \Delta G\|_1 + \|\Delta F \cdot \Delta G\|_1 \\
&\leq \|\Delta F\|_1 \cdot \|G\|_1 + \|F\|_1 \cdot \|\Delta G\|_1 + \|\Delta F\|_1 \cdot \|\Delta G\|_1 \\
&\leq (n+1)^2 2^{s+b} + (n+1)^2 2^{s+b} + (n+1)^2 \\
&\leq (n+1)^2 \cdot 2^{s+b+2}
\end{aligned}$$

holds. For $\tilde{h} := 2^{-2s} \tilde{F}\tilde{G}$, it follows that

$$\|\tilde{h} - h\|_1 \leq 2^{-2s}(n+1)^2 \cdot 2^{s+b+2} \leq 2^{b+2\log(n+1)+2-s} \leq 2^{-\ell},$$

hence an $\ell$-bit-approximation as required can be recovered from the exact product of $\tilde{F}$ and $\tilde{G}$ by mere bitshifts. Since $\|\tilde{F}\|_\infty, \|\tilde{G}\|_\infty \leq 2^{s+b}$, multiplication of $\tilde{F}$ and $\tilde{G}$ can be carried out exactly in $O(\mu((s+b)n))$ bit operations. This proves the complexity result. For the precision requirement, notice that $\|F\|_\infty$, $\|G\|_\infty \leq 2^{s+b}$, and thus we need $(s + b + \lceil\log(n+1)\rceil + 3)$-bit $\|\cdot\|_\infty$-approximations of $f$ and $g$ to compute $\tilde{F}$ and $\tilde{G}$. □

## A.2. Quadratic Interval Refinement for Roots of a Polynomial

Polynomial evaluation is the key operation in many algorithms to approximate the real roots of a square-free polynomial $F(x) \in \mathbb{Z}[x]$: Given an isolating interval $I = (a, b)$ for a real root $\xi$ of $F$ (i.e., $I$ contains $\xi$ and $\bar{I} = [a, b]$ contains no other root of $F$) and an arbitrary positive integer $L$, we aim to compute an approximation of $\xi$ to $L$ bits after the binary point (or, in other words, an $L$-bit approximation of $\xi$) by means of refining $I$ to a width of $2^{-L}$ or less.

A very simple method to achieve this goal is to perform a binary search for the root $\xi$. That is, in the $j$-th iteration (starting with $I_0 := (a_0, b_0) = (a, b)$ in the 0-th iteration), we split the interval $I_j = (a_j, b_j)$ at its midpoint $m(I_j)$ into two equally sized intervals $I'_j = (a_j, m(I_j))$ and $I''_j = (m(I_j), b_j)$. We then check which of the latter two intervals yields a sign change of $F$ at its endpoints,[6] and define $I_{j+1}$ to be the unique such interval. If $F(m(I_j)) = 0$, we can stop because, in this special case, we have exactly computed the root $\xi$. The main drawback of this simple approach is that only linear convergence can be achieved.

In [1], Abbott introduced a method, denoted quadratic interval refinement (QIR), to overcome this issue. It is a trial and error approach which combines the bisection method and the secant method. More precisely, in each iteration, an additional integer $N_j$ is stored (starting with $N_0 = 4$) and (only conceptually) the interval $I_j$ is subdivided into $N_j$ equally sized subintervals $I_{j,1}, \ldots, I_{j,N_j}$. The graph of $f$ restricted to $I_j$ is approximated by the secant $S$ passing through the points $(a_j, F(a_j))$ and $(b_j, F(b_j))$. The idea is that, for $I_j$ small enough, the intersection point $x_S$ of $S$ and the real axis is a considerably good approximation of the root $\xi$, and thus the root $\xi$ is likely to be located in the same of the $N_j$ subintervals as $x_S$. Hence, we compute $x_S$ and consider the unique subinterval $I_{j,\ell}$, with $\ell \in \{1, \ldots, N_j\}$, which contains $x_S$. If $I_{j,\ell}$ yields a sign change of $F$ at its endpoints, we know that it contains $\xi$ and, thus, proceed with $I_{j+1} := I_{j,\ell}$. In addition, we set $N_{j+1} := N_j^2$. This is called a *successful* QIR step. If we are not successful (i.e., there is no sign change of $F$ at the endpoints of $I_{j,\ell}$), we perform a bisection step as above and set $N_{j+1} := \min(4, \sqrt{N_j})$. It has been shown [8] that the QIR method eventually achieves quadratic convergence; in particular, all steps are eventually successful. As a consequence, the bit complexity for computing an $L$-bit approximation of $\xi$ drops from $\tilde{O}(n^3 L)$ (using the bisection approach) to $\tilde{O}(n^2 L)$ (for the QIR method) if $L$ is dominating. Namely, the number of refinement steps reduces from $O(L)$ to $O(\log L)$, and the bit complexity in each step is bounded by $\tilde{O}(n^2 L)$ for both methods (exact polynomial evaluation at a rational number of bitsize $L$).

In [9, 10], a variant of the QIR method, denoted AQIR, has been proposed. It is almost identical to the original QIR method; however, for the sign evaluations and the computation of $x_S$, exact polynomial arithmetic over the rationals has been replaced by approximate but certified interval arithmetic. AQIR improves upon QIR with regard to two main aspects: First, it works for arbitrary real polynomials whose coefficients can only be approximated. Second, it allows to run the computations with an almost optimal precision in each step which is due to an adaptive precision management and the fact that the evaluation points are chosen "away from" the roots of $p$; see [9, 10] for details. In particular, the precision requirement in the worst case drops from $O(nL)$ to $O(L)$ in each step, thus resulting in an overall improvement from $\tilde{O}(n^2 L)$ to $\tilde{O}(nL)$ with respect to bit complexity. Now, if isolating intervals for *all* real roots of $p$ are given, then computing $L$-bit approximations of all real roots uses $\tilde{O}(n^2 L)$ bit operations since we have to consider the cost for the refinement of each of the isolating intervals as many times as the number of real roots (which is at most $n$). *This is the point, where approximate multipoint evaluation comes into play.* Namely, instead of considering the evaluations of $f$ for each interval independently, we can perform $n$ many of these evaluations in parallel without paying more than a polylogarithmic factor compared to only one evaluation. This yields a total bit complexity of $\tilde{O}(nL)$ for computing $L$-bit approximations of all real roots.

---

[6] Here, it is important that $f$ is considered to be square-free. Thus, $\xi$ must be a simple root, and any isolating interval $I = (a, b)$ for $\xi$ yields $F(a) \cdot F(b) < 0$.

We remark that the latter bound is optimal up to logarithmic factors because reading the output already needs $\Theta(nL)$ bit operations.

For the special case, where $p$ has integer coefficients, we recall the following result from Section 3.3:

**Theorem 13.** *Let $F \in \mathbb{Z}[x]$ be a square-free polynomial of degree $n$ with integer coefficients bounded by $2^\tau$, and let $L$ be an arbitrary given positive integer. Then, computing isolating intervals for all real roots of $F$ of width $2^{-L}$ or less uses $\tilde{O}(n^3\tau + nL)$ bit operations.*

*Proof.* Let $\xi_1, \ldots, \xi_m$ denote the real roots of $F$. We proceed in three steps:

In the first step, we compute isolating intervals $I_{\xi_1}, \ldots, I_{\xi_m}$ for all real roots.

In the second step, the intervals are refined such that

$$w(I_{\xi_k}) < w_{\xi_k} := \frac{|F'(\xi_k)|}{32\, e\, d^3 2^\tau \max\{1, |\xi_k|\}^{d-1}} \quad \text{for all } k = 1, \ldots, m, \tag{13}$$

where $e \approx 2.71\ldots$ denotes the Eulerian number. For the latter two steps, we use an asymptotically fast real root isolation algorithm, called NewDsc, which has been introduced in [17]. The proof of [17, Theorem 10] shows that we need $\tilde{O}(n^3\tau)$ bit operations to carry out all necessary computations.

Finally, we use AQIR to refine the intervals $I_{\xi_k}$ to a size of $2^{-L}$ or less. Since the intervals $I_{\xi_k}$ fulfill the inequality (13), [9, Corollary 14] yields that each AQIR-step will be *successful* if we start with $I_0 := I_{\xi_k}$ and $N_0 := 4$. That is, in each of the subsequent refinement steps, $I_j$ will be replaced by an interval $I_{j+1}$ of width $w(I_j)/N_j$, and we have $N_{j+1} = N_j^2$. In other words, we have quadratic convergence right from the beginning and never fall back to bisection. According to [10, Lemma 21] and the preceding discussion, the needed precision for each polynomial evaluation in the refinement steps is bounded by $\tilde{O}(L + n\Gamma_F + \Sigma_F)$, where $2^{\Gamma_F}$ denotes a bound on the modulus of all complex roots $z_1, \ldots, z_n$ of $F$, $\Sigma_F := \sum_{i=1}^n \log \sigma(z_i)^{-1}$, and $\sigma(z_i) := \min_{j \neq i} |z_i - z_j|$ the separation of $z_i$. For a polynomial $F$ with integer coefficients of absolute value $2^\tau$ or less, we may consider $\Gamma_F = 2^{\tau+1}$ according to Cauchy's root bound, and, in addition, it holds that $\Sigma_F = \tilde{O}(n\tau)$; see [10] and the references therein for details. Thus, the bound on the needed precision simplifies to $\tilde{O}(L + n\tau)$. In each iteration of the refinement of a single interval $I_{\xi_k}$, we have to perform a constant number of polynomial evaluations,[7] hence there are $O(n)$ many evaluations for all intervals. All of the involved evaluation points are located in the union of the intervals $I_{\xi_k}$, and thus they have absolute value bounded by $2^\tau$. In addition, $p$ has coefficients of absolute value bounded by $2^\tau$. Hence, in each iteration, we need $\tilde{O}(n^2\tau + nL)$ bit operations for all evaluations according to Theorem 9. Since we have quadratic convergence for all intervals, there are only $O(\log L)$ iterations for each interval, hence the claimed bound follows. □

---

[7] In fact, there are up to 9 evaluations in each step. See [9, Algorithm 3] for details.